\documentclass[12 pt]{article}
\usepackage{EDE}
\allowdisplaybreaks[4]

\renewcommand{\SSS}{\mathbb{S}}
\newcommand{\mbeta}{\boldsymbol{\beta}}

\newcommand{\mPnp}{\mP_{n,\bp}}

\newcommand{\mQnp}{\mu_{n,\bp}}
\newcommand{\pl}{\mP_0(\{y\})}
\usepackage{stmaryrd}

\newcommand{\mY}{\mathbf{Y}}
\newcommand{\my}{\mathbf{y}}

\newcommand{\mf}{\mathbf{f}}

\title{On the bracketing entropy condition and generalized empirical measures}
\author{Davit Varron\\ Université de Bourgogne-Franche-Comté, UMR CNRS 6623}
\begin{document}
\date{}
\maketitle
\textbf{Abstract}: 
 \small{We prove a Donsker and a Glivenko--Cantelli theorem for sequences of random discrete measures generalizing empirical measures. Those two results hold under standard conditions upon bracketing numbers of the indexing class of functions. As a byproduct, we derive a posterior consistency and a Bernstein--von Mises theorem for the Dirichlet process prior, under the topology of total variation, when the observation space is countable. We also obtain new information about the Durst--Dudley--Borisov theorem}.
\section{Introduction}
In this article we shall adopt the generic notation (for $r\in [1,\infty])$
\begin{align*}
\ell^r:=&\aoo \bp\in \RRR^{\NNN},\; \mmi \bp\mmi_r<\infty\aff\text{, where}\\
\mmi \bp\mmi_r^r:=&\sli_{i\in \NNN} \mid p_i\mid^r,\;\text{for } 1\le r<\infty,\text{ and where }\\
\mmi \bp \mmi_{\infty}:=&\sup_{i\in \NNN} \mid p_i\mid\text{, writing }\bp=(p_i)_{i\in \NNN}.
\end{align*}
Let $\po \mfX,\AAA_{\mfX}\pf$ be a measurable space. We shall use the notation
\beq Q(f):=\ili_{\mfX} f dQ\label{notation generique de l'integrale Q(f)}\eeq
for a given signed measure $Q$ in $\AAA_{\mfX}$ with finite total variation, and for each $f\in L^1(Q)$.
Any $\mfX$-valued sequence $\my=(y_i)_{i\in \NNN}$ combined with an element $\bp\in \ell^1$ defines a signed discrete measure on $\po \mathfrak{X},\AAA_{\mfX}\pf$ - with finite total variation - through the following formula:
\beq P_{\my,\bp}:=\sli_{i\in \NNN}p_i\dd_{y_i}.\label{Prbpmy}\eeq
Now substitute $\my$ by a $\AAA_{\mfX}^{\otimes \NNN}$ measurable sequence $\mY=(Y_i)_{i\in\NNN}$, and $\bp$ by a $\ell^1$-valued Borel random variable $\mbeta=(\beta_i)_{i\in \NNN}$ (both of them on a probability space $\wap$). Then the composition map $P_{\mY,\mbeta}$ defines random signed measure in the following sense: for any specified bounded Borel function $f$,  the map
 \beq P_{\mY,\mbeta}(f):\;\;\;\omega\rar  P_{\mY_n(\omega),\mbeta_n(\omega)}(f)\label{Prn}\eeq
 is Borel from $\po \Omega, \AAA\pf$ to $\RRR$. In the sequel we shall continue to adopt the same convention (\ref{notation generique de l'integrale Q(f)}) for $P(f)$ when $P$ is a random or non random measure, and we shall extend it - when meaningful - to functions $f$ that are not necessarily bounded.\lb 
In \cite{Varron14Donsker}, Varron started the investigation on how well known results in empirical processes theory (see, e.g., \cite{Dudley,Vander} for monographs on the subject) could be carried over sequences of random signed measures of the form $P_{\mY_n,\mbeta_n}$
where, for each $n$, the sequence $(Y_{i,n})_{i \in \NNN}$ is independent and identically distributed \textit{given} $\mbeta_n$.
He showed that the uniform entropy numbers and the Koltchinskii--Pollard uniform entropy integral - two crucial notions in empirical processes theory - both adapt very well to that wider class of random measures, which not only encompasses the empirical measure, but also discrete nonparametric Bayesian priors. The latter notion of uniform entropy integral can be briefly defined as follows for a class $\FF$ of real Borel functions on $\po \mathfrak{X},\AAA_{\mfX}\pf$:
\begin{align} \nono J(\dd,\FF):=&\ili_{0}^{\dd}\sqrt{\log\poo \sup_{Q\;probab.} N\po \e\mmi F\mmi_{Q,2},\FF,\norm_{Q,2}\pf\pff}d\e,\;\; \dd\in (0,\infty].
\end{align}
Here $\norm_{Q,2}$ stands for the $L^2(Q)$ norm, $N(\e,\FF,\norm_{Q,2})$ denotes the minimal number of $\norm_{Q,2}$ balls with radius $\e$ needed to cover $\FF$, and $F$ stands for the minimal measurable envelope of the class $\FF$ - see, e.g., \cite[p. 85]{Vander}. Note that $F$ can be simply taken as
$$F(y):=\sup\aoo\mid f(y)\mid, f\in \FF\aff,\; y\in \mfX,$$
when $\FF$ is countable or pointwise measurable - see \S\ref{subsection definition de G_n} below. 
When $J(\infty,\FF)$ is finite, Varron proved a Donsker theorem under natural asymptotic conditions upon $(\mbeta_n,\mY_n)_{n\geq 1}$. Those two asymptotic theorems (see \cite[Theorems 1 and 2]{Varron14Donsker}) involve processes of the form 
\beq G_n(f):=\sli_{i\in \NNN} \beta_{i,n}\cooo f(Y_{i,n})-\EEE\poo f(Y_{i,n})\mid \beta_{i,n}\pff\cfff,\; f\in \FF,\label{definition de Gn}\eeq
indexed by a class $\FF$ of real Borel functions. A rigorous definition of $G_n(\cdot)$ is not immediate and is therefore voluntarily postponed to \S \ref{subsection definition de G_n}.\lb
While the uniform entropy has been celebrated as a very useful condition to prove that a class $\FF$ is Donsker or Glivenko--Cantelli, another condition turned out to be very fruitful as well: bracketing entropy. The bracket $\llbracket f^-,f^+\rrbracket$ between two Borel functions $f^-$ and $f^+$ is defined as the set of Borel functions $f$ fulfilling $f^-\prec f\prec f^+$, the symbol $\prec$ standing for the everywhere pointwise comparison between real functions on $\mfX$.
Denoting by $N_{[]}(\e,\FF,\norm_{Q,2})$ the minimal number of brackets with $\norm_{Q,2}$ diameter less than $\e$ needed to cover $\FF$, the $Q$ bracketing entropy of $\FF$ is defined as 
\beq J_{[]}(\dd,\FF,Q):=\ili_{0}^{\dd}\sqrt{\log N_{[]}(\e,\FF,\norm_{Q,2})}d\e,\;\; \dd\in (0,\infty]\label{definition entropy crochet}.\eeq
A naturally arising question is then: \textit{does bracketing entropy adapt with the same efficiency to sequences of random measures such as in (\ref{Prn})?} The answer provided in the present article is: \textit{yes, but to a lesser extent}. More restrictions upon the weights are needed. First the $\beta_{i,n}$ have to be non negative, since the idea of bracketing relies on the comparison principle
$$f^-\prec f \prec f^+\Rightarrow Q(f^-)\le Q(f)\le Q(f^+), $$ 
when $Q$ is a non negative measure. Second, when looking for a Donsker theorem, $\mmi \mbeta_n\mmi_{\infty}$ has to tend to zero fast enough to counterbalance a the possible growth of $\mmi \mbeta_n\mmi_1$. The amount of compensation is directly linked to the moments of $F(Y_{1,n}),\; n\in \NNN^*$.\lb
Those two conditions were not required under the assumption that $J(\infty,\FF)$ is finite (see \cite[Theorems 1 and 2]{Varron14Donsker}). This difference can be explained by the fact that the use of the Koltchinskii--Pollard entropy is intimately linked to that of symmetrization, namely the study of 
$$G_n^0(f):=\sli_{i\in \NNN}\e_i\beta_{i,n}f(Y_{i,n}),\; f\in \FF,$$
where the $\e_i$ are symmetric Bernoulli (or Rademacher) random variables, independent of $(\mY_n,\mbeta_n)$. By subgaussianity of Rademacher processes, the $G_n^0(\cdot)$ inherit several properties of infinite dimensional Gaussian analysis. In particular, Hilbert spaces take a predominant role. This explains why the results in \cite{Varron14Donsker} hold under conditions upon $\mmi\mbeta_n\mmi_2$ and $\mmi\mbeta_n\mmi_4$. On the other hand, bracketing methods do not rely on subgaussianity, but on a form of Bernstein's inequality. The latter is a tradeoff between subgaussian and subexponential tails for sums of independent random variables that are uniformly bounded. This roughly explains why $\mmi \mbeta_n\mmi_{\infty}$ - and its conjugate norm $\mmi \mbeta_n\mmi_1$ - needs to be controlled. Such a difference of extent between bracketing and uniform entropy was not visible on the empirical process for the following simple reason: when taking $\beta_{i,n}\equiv n^{-1/2}$ for $i\le n$ and $\beta_{i,n}\equiv 0$ otherwise, one has $\mmi \mbeta_n\mmi_{\infty}\equiv\mmi \mbeta_n\mmi_1^{-1}=n^{-1/2}$. This equality makes the counterbalance between those two norms hardly visible in the proof of the bracketing Donsker theorem.\lb
Various interesting classes admit a finite bracketing entropy - see, e.g., \cite[Chapter 2.7]{Vander}. In addition, several examples of posterior distributions in (discrete) Bayesian nonparametrics have the form $\mP_{\mY_n,\mbeta_n}$, or at least exhibit a predominant term that can be expressed as such - see \cite[Section 3]{Varron14Donsker}. Hence our main results present an interesting range of applications, which we will here illustrate through two examples. The first one takes place in the framework of frequentist asymptotic analysis of nonparametric Bayesian priors: for a countable observation space, we prove a posterior consistency and a Bernstein--von Mises theorem for the Dirichlet process prior, under the topology of total variation (see \S \ref{subsection posterior analysis}). Along the proof, we also revisit the Durst--Dudley--Borisov theorem and we obtain additional information about this phenomenon. Our second example of application is a Donsker theorem - under a bracketing condition - for a specific form of local empirical measures (see \S \ref{subsection mesure empirique locale}).
The remainder of this article is organized as follows: in \S \ref{subsection definitions} we give a careful description of the mathematical framework. Then our two main results are stated in \S \ref{section: results}. Applications follow in \S \ref{section: applications}. The proofs of those results are then written in \S \ref{section proofs}. Finally, the Appendix is dedicated to a minor proof.
\section{The mathematical framework}\label{subsection definitions}
In order to properly state our main results, we first need to carefully define their underlying probabilistic framework. This section may be skipped at first reading.
\subsection{The underlying probability space}\label{subsection: definition de wap}
Empirical processes carry over some lacks of measurability that are usually tackled by using outer expectations - see, e.g., \cite[Chapter 1.2]{Vander}.  In order to make use of Fubini's theorem - which is in general untrue for outer expectations - mathematical rigor imposes to define the underlying probability space $\wap$ as a suitable product space. First, for fixed $n\geq 1$, consider a Markov transition kernel from $\ell^1$ to $\mfX$, i.e., a family $\{\mathbf{P}_{n,\bp},\;\bp\in \ell^1\}$ of probability measures for which the maps $\bp\rar \mP_{n,\bp}(A),\; A\in \AAA_{\mfX}$ are measurable from $\po \ell^1,Bor(\ell^1)\pf$ to $\po [0,1],Bor([0,1])\pf$. Also consider a probability measure $Q_n$ on $\po \ell^1,Bor(\ell^1)\pf$ and define:
\begin{align*}
&\tilde{\Omega}:= \ell^1\times\mfX^\NNN,\text{ endowed with its product }\sig\text{-algebra}\\
&\tilde{\AAA}:=Bor(\ell^1)\otimes\AAA_{\mfX}^{\otimes \NNN}
\text{, with probability law defined through the generic formula: }\\
&\PPP_n\poo \ao(\bp,\my)\in \tilde{\Omega}, \bp \in A, \forall j\in \{ 1,\ldots, k\} , y_i\in B_j\af\pff:= \ili_{\bp \in A}\prolijk \mP_{n,\bp}(B_j) dQ_n(\bp).
\end{align*}
Then define $\Omega:=\tilde{\Omega}^{\NNN^*}$, $\AAA:=\tilde{\AAA}^{\NNN^*}$, 
$\PPP:=\bigotimes_{n\geq 1} \PPP_n$ on $\AAA$ and define the $\mY_n$ and $\mbeta_n$ as coordinate maps on $\Omega$~:
\begin{align*}
\mbeta_n(\bp_1,\my_1,\bp_2,\my_2,\ldots):=\bp_n,\text{ and }
\mY_n(\bp_1,\my_1,\bp_2,\my_2,\ldots):=\my_n.
\end{align*}
Note that, for fixed $n$ and $\bp\in \ell^1$, $\mP_{n,\bp}^{\otimes \NNN}$ is the law of $\mY_n$ given $\mbeta_n=\bp$. We shall denote by $\mP_n$ the law of $Y_{1,n}$.\lb
To simplify the notations we now adopt the following convention: each time a map $\mathfrak{h}$ is defined on a probability space, the symbol $\EEE^*(\mathfrak{h})$ will denote the outer expectation with respect to that probability space. We shall adopt the same convention for outer probabilities $\PPP^*$.
\subsection{Definition of $G_n$}\label{subsection definition de G_n}
From now on, and throughout all this article, we shall make the assumption that $\mP_n(F)<\infty$ for all $n\geq 1$.
We also assume that $\FF$ is \textit{pointwise measurable} with countable separant $\FF_0$ in the following sense: for any $f\in \FF$, there exists $(f_m)_{m\geq 1}\in \FF_0^{\NNN^*}$ such that $f_m(y)\rar f(y)$ for each $y\in \mfX$. Such a very standard assumption will be useful to tackle annoying measurability issues. \lb 
Because the symbol $\Sig_{i\in \NNN}$ in (\ref{definition de Gn}) is ambiguous, we need to give a rigorous definition of the processes that will be involved in this article. Our definition differs from that used in \cite{Varron14Donsker} for two reasons. The first (minor) one is to cover the case where the $\mmi \mbeta_n\mmi_1$ are not deterministically equal to 1. The second one is for technical purposes: in our proofs, we shall truncate the $f\in \FF$ from above using thresholds that depend upon the weights. 
First note that, for any bounded function $f$ and any Borel map $T$ from $\ell^1$ to $\RRR^+$, the map
$$\Phi_{f,T}:\; (\my,\bp)\rar \sli_{i\in \NNN}p_i\poo f\indic_{\{F\le T(\bp)\}}(y_i)-\mP_{n,\bp}\po f\indic_{\{F\le T(\bp)\}}(y_i)\pf\pff$$
is properly defined (through the limits in $\RRR$ of partial sums) and Borel from $\mfX^{\NNN}\times \ell^1$ to $\RRR$. One can hence define a random variable $G_n^T(f)$ by composition $G_n^T(f):=\Phi_{f,T}\circ(\mY_n,\mbeta_n)$. 
We will say that a map $\psi$ from $\FF$ to $\RRR$ is \textit{$\FF_0$-separable} whenever we have $\mmi \psi\mmi_{\FF}=\mmi \psi \mmi_{\FF_0}$. We shall also denote by $\BB(\FF,\FF_0)$ the space of all bounded $\FF_0$-separable functions and by $\AAA_{\norm_{\FF}}$ is the $\sig$-algebra spanned by the $\norm_{\FF}$-balls.
\begin{lem}
For any choice of $T$ as above, the map
$$\Phi_{\FF,T}:\; (\my,\bp)\rar \ao f\rar \Phi_{f,T}(\my,\bp)\af$$
is measurable from $\mfX^{\NNN}\times \ell^1$ to $(\BB(\FF,\FF_0),\AAA_{\norm_{\FF}})$.
\end{lem}
\textbf{Proof}:
Fix $T$. Let us first prove that $\Phi_{\FF,T}$ takes its values in $\BB(\FF,\FF_0)$. Fix $\my\in \mfX^{\otimes \NNN}$ and $\bp\in \ell^1$. Since, for all $k\in \NNN$:
$$\sup_{f\in \FF}\Mid\sli_{i\geq k+1}p_i \coo f\indic_{\{F\le T(\bp)\}}(y_i)-\mP_{n,\bp}\po f\indic_{\{F\le T(\bp)\}}\pf\cff\Mid \le 2T(\bp)\sli_{i\geq k+1}\mid p_i\mid,$$
and since $\po \BB(\FF,\FF_0),\norm_{\FF}\pf$ is a Banach space, it is sufficient to prove that each trajectory $f\rar f\indic_{\{F\le T(\bp)\}}(y_i),\;i\in \NNN,$ and $f\rar \mP_{n,\bp}\po f\indic_{\{F\le T(\bp)\}}\pf$ is $\FF_0$-separable. To see this, take $f\in \FF$ and consider $(f_m)_{m\geq 1}\in \FF_0^{\NNN}$ such that $f_m\rar f_0$ pointwise. Thus $\mP_{n,\bp}(f_m\indic_{\{F\le T(\bp)\}})\rar \mP_{n,\bp}\po f\indic_{\{F\le T(\bp)\}}\pf$ by the dominated convergence theorem. Now since $\Phi_{\FF,T}$ takes its values in $\BB(\FF,\FF_0)$, the following equality holds for any $(\my,\bp)\in\mfX^{\NNN}\times \ell^1$:
\beq \sup_{f\in \FF}\mid \Phi_{f,T}(\my,\bp)\mid=\sup_{f\in \FF_0}\mid \Phi_{f,T}(\my,\bp)\mid,\eeq
and then the measurability of each $\Phi_{f,T},\; f\in \FF_0$ ensures that of $\Phi^T_{\FF}$ with respect to $\AAA_{\norm_{\FF}}$. $\Box$\lb
Now denote by $\tilde{\EE}_{\FF,\FF_0}$ the space of all measurable maps from $\po \Omega,\AAA\pf$ to $\po \BB(\FF,\FF_0),\AAA_{\norm_{\FF}}\pf$.
The preceding lemma gives the opportunity to define the processes $G_n$ on any class of functions $$\FF_M:=\ao f\indic_{\{F\le M\}},\;f\in \FF\af$$ since (confounding $M$ with a constant function on $\ell^1$) the composition map
$$G_n^M:=\Phi_{\FF,M}\circ(\mY_n,\mbeta_n)$$
belongs to $\tilde{\EE}_{\FF,\FF_0}$. Denote by $\EE_{\FF,\FF_0}$ the quotient space of $\tilde{\EE}_{\FF,\FF_0}$ with respect to the equivalence class 
$$G\sim G'\Leftrightarrow \mmi G-G'\mmi_{\FF}=0,\; \PPP\text{-a.s.} ,$$
and endow $\EE_{\FF,\FF_0}$ with the compatible distance 
\beq d(G,G'):=\EEE\poo \arctan \po \mmi G-G'\mmi_{\FF}\pf\pff,\label{definition de la distance L_0 processus}\eeq
which is that of $\norm_{\FF}$-convergence in probability.
The following lemma defines $G_n$ as a suitable limit of the $G_n^M$ when $M\rar \infty$.
\begin{lem}\label{lem: definition complete de Gn}
For fixed $n\geq 1$, the sequence $(G_n^M)_{M\geq 1}$ is Cauchy in the complete metric space $\po\EE_{\FF,\FF_0},d\pf$. It hence converges to a limit which we take as the definition of $G_n$. Moreover, for any sequence $(T_k)$ of Borel thresholding maps fulfilling $T_k(\mbeta_n)\cvproba \infty$ as $\kif$, we have 
$d(G_n^{T_k},G_n)\rar 0$ as $\kif$.
\end{lem}
\textbf{Proof}: 
For integers $M,M'$ we have, writing $f^{M,M'}:=f\indic_{\{M<F\le M'\}}$
\begin{align}
\nono  &d\poo G_n^M,G_n^{M'}\pff \\
\nono = &\EEE\poooo \arctan \pooo \sup_{f\in \FF}\Mid \Phi_{f^{M,M'}}(\mY_n,\mbeta_n)\Mid \pfff\pffff   \\
 \le & \EEE\poooo \arctan \poooo \sli_{i\in \NNN}\mid \beta_{i,n}\mid\cooo F^{M,M'}(Y_{i,n})+\EEE\poo F^{M,M'}(Y_{i,n})\Mid \mbeta_n\pff\cfff\pffff\pffff.\label{pnop}
 \end{align}
Using Fatou's lemma for conditional expectations and the concavity of $\arctan$ on $\RRR^+$ we have, almost surely:
\begin{align}
\nono&\EEE\poooo \arctan \pooo \sli_{i\in \NNN}\mid \beta_{i,n}\mid\cooo F^{M,M'}(Y_{i,n})+\EEE\poo F^{M,M'}(Y_{i,n})\Mid \mbeta_n\pff\cfff\pfff\Mid \mbeta_n\pffff\\
\nono\le&\arctan\pooo2\sli_{i\in \NNN} \mid \beta_{i,n}\mid \EEE\poo F^{M,M'}(Y_{i,n})\mid \mbeta_n\pff\pfff\\
\label{jil}=&\arctan\pooo2\mmi \mbeta_{n}\mmi_1\EEE\poo F^{M,M'}(Y_{1,n})\mid \mbeta_n\pff\pfff\\
\le &\arctan\pooo2\mmi \mbeta_{n}\mmi_1\EEE\poo F\indic_{\{F>M\}}(Y_{1,n})\mid \mbeta_n\pff\pfff,\label{bepplap}
\end{align}
where $(\ref{jil})$ comes from the fact that the law $\mY_n$ given $\mbeta_n=\bp$ is $\mP_{n,\bp}^{\otimes \NNN}$.
It hence suffices to prove that the right hand side (RHS) of (\ref{bepplap}) tends to $0$ in probability as $M\rar \infty$. This is true since $\mP_n(F)<\infty$. Now to prove the last statement of Lemma \ref{lem: definition complete de Gn}, formally replace $M$ by $T_k(\bp)$ in the preceding calculus and let $M'\rar\infty$ to obtain, using Fatou's lemma for conditional expectations:
\beq
d\po G_n^{T_k},G_n\pf
\le \EEE\poooo \arctan\pooo 2\mmi \mbeta_{n}\mmi_1\EEE\poo F\indic_{\{F>T_k(\mbeta_n)\}}(Y_{1,n})\mid \mbeta_n\pff\pfff\pffff\label{inegalite entre les esperance condtionnelles pour le reste de cauchy},\eeq
which tends to $0$ as $\kif$, by assumption upon $(T_k(\mbeta_n))_{k\geq 1}$ and since $\mP_n(F)<\infty$.
 $\Box$
 \subsection{Definition of the limit processes}
One of our results is a Donsker theorem (see Theorem \ref{T2} below). The limit processes are mixtures of $\FF$-indexed Brownian bridges, for which a rigorous definition is not immediate due to the non separability of $\ell^{\infty}(\FF)$. We shall use the definition of \cite{Varron14Donsker}. However, since the possibility to condition upon the weights was not made perfectly clear in \cite{Varron14Donsker}, we feel the need to give more precisions in its reminder. First fix $n\geq 1$.
For $p\geq 1$, $\mf:=(f_1,\ldots,f_p)\in \FF^p$ and $\bp\in\ell^1$, write $\mQnp^\mf$ for the centered Gaussian distribution on $\RRR^p$ with covariance matrix 
$$\Sig_{n,\bp}^\mf:=\coo \mPnp\poo (f_j-\mPnp(f_j))(f_{j'}-\mPnp(f_{j'}))\pff\cff_{(j,j')\in \{1,\ldots, p\}^2}.$$ 
Now consider $\Omega':=\RRR^{\FF}$, endowed with its product Borel $\sig$-algebra $\AAA'$, i.e the $\sig$-algebra spanned by the $\pi$-system 
\begin{align}
\nono &\AAA'_0:=\aoo C_{A_{f_1},\ldots,A_{f_k}},\: k\geq 1,\; \mf\in \FF^k,\; A_{f_j}\text{ Borel for each }j\in \{1,\ldots,k\}\aff\text{, with}\\
\nono &C_{A_{f_1},\ldots,A_{f_k}}:=\aoo \psi \in \Omega',\; \forall j\in \{1,\ldots,k\},\; \psi(f_j)\in A_{f_j}\aff.
\end{align}
For fixed $\bp\in \ell^1$, Kolmogorov's extension theorem (see, e.g., \cite[p. 115, Theorem 6.16]{Kallenberg}) ensures the existence of a unique probability measure $\PPP'_{n,\bp}$ on $\po \Omega',\AAA'\pf$ which is compatible with the system of marginals $\ao\mQnp^{\mf},\;\mf\in\FF^p,\; p\geq 1\af$, namely $$\PPP'_{n,\bp}\po C_{A_{f_1},\ldots,A_{f_k}}\pf=\mQnp^\mf\po A_{f_1}\times\ldots\times A_{f_k}\pf,$$
for all elements of $\AAA'_0$. Now since $\{\mP_{n,\bp},\;\bp\in \ell^1\}$ defines a transition kernel, then so does the family $\{\mQnp^\mf,\;\bp \in \ell^1\}$ for any specified $\mf$: it is a simple consequence of the fact that $\bp \rar \Sig_{n,\bp}^{\mf}$ is Borel. Consequently the map $\bp \rar \PPP'_{n,\bp}(C)$ is measurable for any specified $C\in \AAA'_0$. That measurability property is then immediately extended to all $C\in \AAA$ by Dynkin's $\pi$-$\lab$ theorem - see, e.g., \cite[p. 2, Theorem 1.1]{Kallenberg}. As a consequence, the family $\{\PPP'_{n,\bp},\;\bp\in \ell^1\}$ defines a transition kernel on $\po \Omega',\AAA'\pf$. Finally, we take $\Omega'':=\ell^1\times \Omega'$, endowed with $Bor(\ell^1)\otimes \AAA'$, and with probability $\PPP''_n$ defined by 
$$\PPP'_n(A\times C):=\ili_{\bp \in A}\PPP'_{n,\bp}(C)dQ_n(\bp).$$
 We then take $\AAA''_n$ as the completion of $Bor(\ell^1)\otimes \AAA'$ with respect to $\PPP''_n$ and we define $W_n$ as the coordinate map $(\bp,\psi)\rar \psi$ on  $\po \Omega'',\AAA''_n,\PPP_n''\pf$. \begin{lem}
Let $n\geq 1$ be an integer. Assume that, $\PPP$-almost surely:
$$J_{[]}\poo\infty,\FF,\norm_{\mP_{n,\mbeta_n},2}\pff<\infty.$$
Then $W_n$ is $\PPP''_n$-almost surely bounded. 
\end{lem}
\textbf{Proof:} Write $\mbeta'_n$ as the canonical map $(\bp,\psi)\rar \bp$ on $\po \Omega'',\AAA''_n,\PPP''_n\pf$ and note that $\mbeta'_n$ and $\mbeta_n$ are equal in law. For any finite subclass $\{f_1,\ldots,f_p\}\subset \FF$ we have, by conditioning on $\mbeta_n'$ and then using Dudley's chaining theorem (see, e.g., \cite[p. 101, Corollary 2.2.8]{Vander})
\begin{align*}
\EEE\pooo \max_{j\le p}\mid W_n(f_j)\mid\Mid \mbeta'_n\pfff\le &\mathfrak{C}_0
\ili_{0}^{+\infty}\sqrt{\log N\po \e,\FF,\norm_{\mP_{n,\mbeta'_n},2}\pf}d\e\\
\le &\mathfrak{C}_0
\ili_{0}^{+\infty}\sqrt{\log N_{[]}\po \e,\FF,\norm_{\mP_{n,\mbeta'_n},2}\pf}d\e\\
=&
\mathfrak{C}_0J_{[]}\poo\infty,\FF,\mP_{n,\mbeta'_n}\pff,\;\PPP''_n\text{-almost surely},
\end{align*}
where $\mathfrak{C}_0$ is a universal constant. It follows by Fatou's lemma for conditional expectations that
$$\EEE\pooo \sup_{f\in \FF_0}\mid W_n(f)\mid \Mid \mbeta'_n\pfff\le \mathfrak{C}_0J_{[]}\poo\infty,\FF,\mP_{n,\mbeta'_n}\pff,\;\PPP''_n\text{-almost surely}.$$
Now using the same arguments as those used to obtain \cite[p. 2314, assertion (49)]{Varron14Donsker}, one can show that
$$\PPP''_n\poo \ao \omega\in \Omega'',\; W_n(\omega)\text{ is not }\FF_0\text{ separable}\af\pff=0,$$
which concludes the proof.$\Box$ \lb

\section{Results}\label{section: results}
Before stating our two main results, let us briefly mention that the maps $\bp\rar N_{[]}\po\e,\FF,\norm_{\mP_{n,\bp},r}\pf$ and $\bp\rar J_{[]}\po \dd,\FF,\mP_{n,\bp}\pf$ are properly measurable for fixed $\e$ and $\dd$. This is proved in \S
\ref{sous section preuve de la mesurabilite}.
\subsection{A Glivenko--Cantelli theorem}
Our first result is a Glivenko--Cantelli theorem. Recall that $\mP_n$ is the law of $Y_{1,n}$. We shall denote by $\ell^{1,+}:=\ell^1\cap [0,\infty[^{\NNN}$ the set of non negative summable sequences.
\begin{theo}\label{T1}
Assume that 
\beq \lim_{M\rar \infty}\; \lsn \mP_n\poo F\indic_{\{F\geq M\}}\pff=0\label{condition envelope GC},\eeq
and that, for any $\e>0$:
\beq    \pooo N_{[]}\po\e,\FF,\norm_{\mP_{n,\mbeta_n},1}\pf\pfff_{n\geq 1}\text{ is bounded in probability}.\label{condition bracket GC}\eeq
Also assume that $\mbeta_n\in \ell^{1,+}$ is almost surely for all $n$, and that  
\beq \po \mmi \mbeta_n\mmi_1\pf\suite\text{ is bounded in probability.}\label{condition GC mbeta_n borne en proba}\eeq Then, under the condition $\mmi \mbeta_n\mmi_2\cvproba 0$ we have 
$\mmi G_n\mmi_{\FF}\cvproba 0.$
\end{theo}

\textbf{Remark}: It is important to compare the assumptions of Theorem \ref{T1} to 
those of Dudley's bracketing Glivenko--Cantelli theorem \cite[p. 122, Theorem 2.4.1]{Vander} for the empirical measure 
$P_{\mY_n,\mbeta^{Emp}_n}$, where $\beta_{i,n}^{Emp}:\equiv n^{-1}$ for $i\le n$ and is identically null otherwise, and where $\mY_n$ is constant in $n$ - hence with constant law $\mP_n=\mP_0$. A few easy arguments then show that, in this special case, those two theorems exactly coincide: first, for $\mbeta_n=\mbeta^{Emp}_n$ the convergence in probability is equivalent to an almost sure convergence by Pollard's reverse martingale argument (see, e.g., \cite[p. 124, Lemma 2.4.5]{Vander}). Second, (\ref{condition envelope GC})+(\ref{condition bracket GC}) is here equivalent to the finiteness of $N_{[]}\po \e,\FF,\norm_{\mP_0,1}\pf$ for each $\e>0$.
\subsection{A Donsker theorem}
For a sequence $Z_n$ of maps from $\Omega$ to $\RRR$ we shall write
$${\lsn}^{\PPP^*}\; Z_n:=\inf\aoo M\in \RRR,\; \limn \PPP^*(Z_n\geq M)=0\aff,$$
with the convention $\inf_{\emptyset}=+\infty$, and we shall simply write $\lsn^{\PPP}\;Z_n$ when the maps $Z_n$ are measurable.
Our second result is a Donsker theorem. \begin{theo}\label{T2}
Assume that
\begin{align}
&\mmi \mbeta_n\mmi_2\cvproba 1 \label{ condition Donsker norme L2 egale a 1},\\
&\mmi\mbeta_n\mmi_{\infty}\cvproba 0,\label{condition Donsker linfty tend vers 0}
\end{align}
and that, for some $p\in [2,\infty[$ 
\begin{align}
&\mmi\mbeta_n\mmi_1 \times \mmi \mbeta_n\mmi_{\infty}^{p-1}\text{ is bounded in probability},\label{condition Donsker produit l1linfty borne}\\
&\lim_{\dd\rar 0}\; {\lsn}^{\PPP} J_{[]}(\dd,\FF,\mP_{n,\mbeta_n})=0,\label{condition Donsker sur l'entropie crochet}\\
&\lim_{M\rar \infty}\; \lsn \mP_n\po F^p\indic_{\{F>M\}}\pf=0\label{ condition envelope Donsker}.
\end{align}
Also assume that there exists a semimetric $\rho$ that makes $\FF$ totally bounded, and fulfilling
\beq \lim_{\dd\rar 0}\;{\lsn}^{\PPP^*}\mathop{\sup_{(f_1-f_2)\in \FF^2,\;}}_{\rho(f_1,f_2)<\dd}\mP_{n,\mbeta_n}\poo (f_1-f_2)^2\pff=0\label{ condition sur la semimetrique rho}.\eeq
Then
\beq d_{BL}\poo G_n,W_{n}\pff:=\sup_{B\in BL1}\Mid \EEE^*\poo B\po G_n\pf\pff-\EEE^{*}\poo B\po W_n\pf\pff\Mid\rar 0,\label{convergence BL1}\eeq
where $BL1$ is the set of all 1-Lipschitz functions on $\po \ell^{\infty}(\FF),\norm_{\FF}\pf$ that are bounded by 1.\lb
Moreover, if $\FF$ is uniformly bounded, then (\ref{convergence BL1}) holds without assuming (\ref{condition Donsker produit l1linfty borne}) nor (\ref{ condition envelope Donsker}).
\end{theo}
\textbf{Remarks}:
We chose to state Theorem \ref{T2} under the most general assumptions that our methodology can afford. In order to give more substance to those conditions, it seems convenient to discuss on the place of Theorem \ref{T2} in the existing literature on Donsker theorems for empirical processes. \begin{enumerate}
\item When $\mbeta_n$ is the vector of rescaled empirical weights ($\beta_{i,n}\equiv n^{-1/2}$ for $i\le n$ and $\beta_{i,n}\equiv0$ otherwise), and when $\mP_n=\mP_0$ is constant in $n$, then $\mP_{\mY_n,\mbeta_n}$ is a sequence of empirical processes. Noting that - for $p=2$ - the $\mbeta_n$ obviously satisfy conditions (\ref{ condition Donsker norme L2 egale a 1}), (\ref{condition Donsker linfty tend vers 0}) and  (\ref{condition Donsker produit l1linfty borne}) one can immediately conclude that - in this setup - Theorem \ref{T2} exactly coincides with Ossiander's bracketing Donsker theorem \cite[Theorem 3.1]{Ossiander87}. Andersen \textit{et al.} \cite{AndersenGOZ88} did also prove a Donsker theorem under more general conditions, where the finiteness of $J_{[]}(\infty,\FF,\mP_0)$ is relaxed a to more abstract assumption, involving majorizing measures on $\norm_{\mP_0,2}$ balls and "weak $\norm_{\mP_0,2}$" brackets. This possible extension of Theorem \ref{T2} is beyond the scope of the present article and may deserve future investigations.
\item Let us now relax the assumption that $\mP_n$ is constant in $n$. In that case the $G_n$ fall into the framework of triangular arrays of empirical processes with varying baseline measures, which were studied by Sheehy and Wellner \cite[Section 3]{SheehyW92}. These authors did prove a Donsker result for $G_n$ indexed by classes fulfilling $J(\infty,\FF)<\infty$, under the envelope condition (\ref{ condition envelope Donsker}), and assuming that $\mP_n$ converges to a limit $\mP_0$ in the following sense - see their Corollary 3.1: 
\begin{align}
\nono &\sup_{(f_1,f_2)\in \FF^2}\max\aoo\mid\mP_n\po(f_1-f_2)^2\pf-\mP_0\po(f_1-f_2)^2\pf\mid,\; \mid\mP_n(f_1)-\mP_0(f_1)\mid,\\
&\;\;\;\;\;\;\;\;\;\;\;\;\;\;\;\;\;\;\;\;\;\;\;\mid\mP_n(f_1^2)-\mP_0(f_1^2)\mid\aff\rar 0\label{SW}.\end{align}
 It is then clear that our Theorem 2 puts forward an analogue of their result, replacing their assumption $J(\infty,\FF)<\infty$ by the bracketing condition (\ref{condition Donsker sur l'entropie crochet}). To see this, just note that (\ref{ condition sur la semimetrique rho}) is satisfied under (\ref{SW}), by choosing $\rho(f_1,f_2):=\mmi f_1-f_2\mmi_{\mP_0,2}$. 
\item Let us now discuss on assumption (\ref{condition Donsker sur l'entropie crochet}), which might be the most cumbersome to verify for applications. If $J_{[]}(\infty,\FF,\mP_0)<\infty$, a simple way to check (\ref{condition Donsker sur l'entropie crochet}) - by direct comparison of bracketing numbers - is to prove that 
\beq\sup_{(f_1,f_2)\in \FF^2}\frac{\mP_{n,\mbeta_n}\po  (f_1-f_2)^2\pf}{\mP_0\po  (f_1-f_2)^2\pf}\label{condition de domination}\text{ is bounded in probability.}\eeq
 Such a sufficient condition is quite restrictive and seems far from necessary, but its verification is sometimes very simple to perform. This is for example the case for our application to the local empirical process at fixed point - see \S \ref{subsection mesure empirique locale}. \item We conclude this series of remarks by pointing out that, whereas involving random weights, Theorem 2 has almost no connections with Donsker theorems for bootstrap empirical measures. For more details see \cite[Remark 2.2]{Varron14Donsker}. 
\end{enumerate}

\section{Applications}\label{section: applications}
\subsection{Posterior analysis of the Dirichlet process prior under the discrete total variation}\label{subsection posterior analysis}
Assume (in this subsection only) that $\mfX$ is infinite countable. The class $\FF$ of all indicator functions of subsets of $\mfX$:
\beq\FF:=\aoo \indic_{C},\: C\subset \mfX \aff\label{definition classe de tous les sous ensembles}\eeq
is rich enough to define the discrete total variation between two measures on $\mfX$, since
$$\mmi Q-Q'\mmi_{\FF}=\sup_{C\subset \mfX}\mid Q(C)-Q'(C)\mid=:\mmi Q-Q'\mmi_{Tot. var.}.$$
Clearly, $\FF$ is too large to satisfy $J(\infty,\FF)<\infty$. It was however shown in the celebrated Durst--Dudley--Borisov theorem that $\FF$ may have a finite bracketing entropy $J_{[]}(\infty,\FF,Q)$ under a simple necessary and sufficient criterion upon $Q$, namely
\begin{tabbing}
 $ \hskip 0pt $ \= $(DDB(Q))\;:$  $ \hskip 20pt $ \= $\sli_{y\in \mfX}\sqrt{Q(\{y\})} <\infty$.
 \end{tabbing}
 \begin{theo}[Durst--Dudley--Borisov, 1981]\label{theo Durst}
For the class $\FF$ defined in (\ref{definition classe de tous les sous ensembles}) we have
$$J_{[]}(\infty,\FF,Q)<\infty\;\;\Leftrightarrow\;\; \po DDB(Q)\pf.$$
\end{theo}
We shall combine Theorem \ref{T1} with a refinement of Theorem \ref{theo Durst} - see Lemma \ref{lem des sommes de queues} - to prove both a posterior consistency and a Bernstein--von Mises theorem for the Dirichlet process prior, under the \textit{discrete} total variation. To properly state it, we need to introduce some more notations. From now on we shall denote by $DP(\alp,M)$ a Dirichlet process with mean probability measure $\alp$ on $\mfX$ and concentration parameter $M>0$. A possible representation of $DP(\alp,M)$ is that of Sethuraman \cite{Sethuraman94}:
\beq DP(\alp,M)=_{law}Pr_{\mY,\mbeta},\label{representation de Sethuraman}\eeq
where $\mY\leadsto \alp^{\otimes \NNN}$, and $\beta_i:=V_i\;\Pi_{j\le i-1}(1-V_j)$
with $(V_i)_{i\in \NNN}\leadsto Beta(1,M)^{\otimes \NNN}$ being independent of $\mY$.
 Now consider the nonparametric Bayesian model where the prior $Pr$ has distribution $DP(\alp,M)$ and where the sample $\XXN$ has conditional law $\mP^{\otimes n}$ given $Pr=\mP$. In this model it is well known (see \cite{Ferguson73}) that a natural expression of the posterior distribution of $Pr$ given $\XXN=\xxn$ is $Post_n\xxn:=DP(\alp_{\xxn},M+n)$, where 
\begin{align}
 \nono \alp_{\xxn}:=&\tht_n\alp+(1-\tht_n)\mP_{\xxn}\text{, with }
\nono  \tht_n:=\frac{M}{M+n},\text{ and } \mP_{\xxn}:=\frac{1}{n}\sliin \dd_{x_i}.
 \end{align}
We shall take advantage of that explicit representation to prove the following two results.
\begin{coro}[Posterior consistency]\label{ coro:posterior}
Take $\FF$ as in (\ref{definition classe de tous les sous ensembles}).
Let $\mP_0$ be a probability measure on countable $\mfX$, and let $(X_n)_{n \geq 1}\leadsto \mP_0^{\NNN^*}$. Then for almost every sequence $(x_n)_{n\geq 1}$ we have 
$$\Mmi  Post_n\xxn - \mP_{\xxn}\Mmi_{Tot. Var.}\cvloi 0.$$
\end{coro}

\begin{coro}[Bernstein--von Mises]\label{coro: posterior et bernstein}
Assume in addition $(DDB(\mP_0))$ and $(DDB(\alp))$. Then for almost every sequence $(x_n)_{n\geq 1}$ we have 
\beq \nono\sqrt{n} \pooo Post_n\xxn - \mP_{\xxn}\pfff\cvloi \GGG_{\mP_0},\text{ in }\ell^{\infty}(\FF).\eeq
 As a consequence, for almost every sequence $(x_n)_{n\geq 1}$ we have 
$$\sqrt{n}\Mmi  Post_n\xxn - \mP_{\xxn}\Mmi_{Tot. Var.}\cvloi \mmi\GGG_{\mP_0}\mmi_{\FF}.$$
\end{coro}
The corresponding proofs are written in \S \ref{sous section preuve Bernstein}.
 \subsection{A Donsker theorem for local empirical measures under  a bracketing condition}\label{subsection mesure empirique locale}
Assume in this subsection that $\mfX=\RRR^d$.
The local empirical process indexed by functions, introduced by Einmahl and Mason \cite{EinmahlM97} has been intensively investigated during the last decades, due to its connections with several smoothing nonparametric methods. One of its particular forms can be written as follows:
$$T_{n,h_n}(f):=\frac{1}{\sqrt{nh_n^d}}\sliin \coooo f\poo h_n^{-1}\po Z_i-z\pf\pff-\EEE\pooo f\poo h_n^{-1}\po Z_i-z\pf\pff\pfff\cffff,$$
where $(h_n)_{n\geq 1}$ is a deterministic non negative sequence tending to $0$, and  where $(Z_n)_{n\geq 1 }$ is an i.i.d. sequence. Implicit in the results of Einmahl and Mason \cite[Theorem 1.1]{EinmahlM97} is the following Donsker theorem.
\begin{theo}[from Einmahl and Mason, 1997]\label{ theo de Einmahl et Mason}
Let $(h_n)\suite$ be a non random sequence of non negative numbers such that $h_n\rar 0$ and $nh_n^d\rar \infty$.
Assume that $J(\infty,\FF)<\infty$. Assume that the support $S$ of $F$ is bounded, and that $Z_1$ admits a version of Lebesgue density $\mf$ on a neighborhood of $z$ that is continuous at $z$ and such that $\mf(z)>0$. Also assume that, taking $\mP_0$ as the uniform distribution on $S$, we have $\mP_0(F^2)<\infty$. Then we have the following weak convergence \beq  \nono \frac{1}{\sqrt{\lab(S)\mf(z)}}T_{n,h_n}(\cdot)\cvloi \WW_{\mP_0}(\cdot),\text{ in }\ell^{\infty}(\FF), \eeq
where  $\WW_{\mP_0}(\cdot)$ denotes the $L^2(\mP_0)$-isonormal Gaussian process indexed by $\FF$ (or $\mP_0$-Brownian motion).
\end{theo}
Their proof heavily relies on a representation of their own \cite[Proposition 3.1]{EinmahlM97}:
$$T_{n,h_n}(f):=_{law}\frac{1}{\sqrt{nh_n^d}}\times\cooo \sliin b_{i,n}\po f(Y_{i,n})-\mP_n(f)\pf\cfff+R_n(f),$$
as processes indexed by $\FF$, where:
\begin{itemize}
\item The $(b_{i,n})_{i\le n}$ are i.i.d Bernoulli with parameter $a_n:=\PPP\po h_n^{-1}(Y_1-z)\in S\pf$;
\item The $(Y_{i,n})_{i\le n}$ are i.i.d with law 
\beq \mP_n:=\PPP\po Y_1\in \cdot\mid h_n^{-1}(Y_1-z)\in S\pf,\label{definition Pn local}\eeq 
with $(b_{i,n})_{i\le n}\indep (Y_{i,n})_{i\le n}$; 
\item The term
\beq R_n(f):= \frac{\sliin (b_{i,n}-a_n)}{\sqrt{nh_n^d}}\mP_n(f)\label{definition de R_n(f)}\eeq
plays the asymptotic role of a correcting drift between the Brownian bridge $\GGG_{\mP_0}$ and the Brownian motion $\WW_{\mP_0}$.\end{itemize}
The following corollary of Theorem \ref{T2} is a Donsker theorem for $T_{n,h_n}$ under the condition $J_{[]}(\infty,\FF,\mP_0)<\infty$. Its proof is written in \S \ref{sous section preuve de Donsker local}.
\begin{coro}\label{coro Donsker local}
Theorem \ref{ theo de Einmahl et Mason} still holds if assumption $J(\infty,\FF)<\infty$ is replaced by $J_{[]}(\infty,\FF,\mP_0)<\infty.$
\end{coro}

\section{Proofs}\label{section proofs}
\subsection{Proof of Theorem \ref{T1}}
The proof is divided in two lemmas.

\begin{lem}\label{ lem: GC sous un seuil arbitraire}
Take $M>0$. Under the assumptions of Theorem \ref{T1} we have
$\mmi G_n^M\mmi_{\FF}\cvproba 0$ as $\nif$.
\end{lem}
\textbf{Proof }: Fix $M>0,\; \e>0$, and choose - by (\ref{condition bracket GC}) and (\ref{condition GC mbeta_n borne en proba}) - an integer $N$ for which, $\PPP(\mbeta_n\in \SSS_n)>1-\e$ for all $n\geq 1$, with
$$\SSS_n:= \aoo \bp \in \ell^{1,+},\;\mmi \bp\mmi_1\le N\text{ and } N_{[]}\po \e,\FF,\norm_{\mP_{n,\bp},1}\pf\le  N\aff.$$
Now fix $n\geq 1$, and $\bp \in \SSS_n$ and denote by $Br(n,\bp)=\ao(f_j^-,f_j^+),\;j=1,\ldots,N\af$ a covering bracket of $\FF$ with $\max_{j\le N}\mmi f_j^+-f_j^-\mmi_{\mP_{n,\bp},1}\le \e/N$. Using the same comparison argument as in \cite[p. 122]{Vander} we have, for $(f^-,f^+)\in Br(n,\bp)$, $f\in \llbracket f^-,f^+\rrbracket$ and $\my\in \mfX^\NNN$ (recall (\ref{Prbpmy})):
\begin{align} 
\nono &\coo P_{\my,\bp}\poo f^-\indic_{\{F\le M\}}-\mP_{n,\bp}\po f^-\indic_{\{F\le M\}}\pf\pff \cff- N\mP_{n,\bp}( f^+ -f^-)\\
\nono \le & P_{\my,\bp}\poo f\indic_{\{F\le M\}}-\mP_{n,\bp}\po f\indic_{\{F\le M\}}\pf\pff\\
 \le&  \coo P_{\my,\bp}\poo f^+\indic_{\{F\le M\}}-\mP_{n,\bp}\po f^+\indic_{\{F\le M\}}\pf\pff\cff + N\mP_{n,\bp}( f^+ -f^-),\label{comparison 1}
\end{align}
from where, writing $B_{n,\bp}$ for the set of functions $f$ that are a side of a bracket (hence $\sharp B_{n,\bp}\le 2N$, where "$\sharp$" stands for "cardinal"):
\begin{align}
\nono &\sup_{f\in \FF}\Mid P_{\my,\bp}\poo f\indic_{\{F\le M\}}-\mP_{n,\bp}\po f\indic_{\{F\le M\}}\pf\pff\Mid\\
\nono \le&  \max_{f\in B_{n,\bp}}\Mid P_{\my,\bp}\poo f\indic_{\{F\le M\}}-\mP_{n,\bp}\po f\indic_{\{F\le M\}}\pf\pff\Mid+\e.\end{align}
Note that the condition $\bp\in \ell^{1,+}$ is crucial to obtain (\ref{comparison 1}).
Now formally replacing $\my$ by an i.i.d sequence $(Y_i)_{i\in \NNN}$ having distribution $\mP_{n,\bp}$ we obtain, for $\bp\in\SSS_n$:
\beq \EEE\poooo \sup_{f\in \FF}\Mid \sli_{i\in \NNN}p_i\po f(Y_i)-\mP_{n,\bp}(f)\pf\Mid \pffff\le \Delta_n(\bp)+\e,\text{ where }\label{ha}\eeq
\begin{align}
\nono\Delta^2_n(\bp):=&\EEE^2\poooo \sli_{f\in B_{n,\bp}}\Mid \sli_{i\in \NNN}p_i\poo f(Y_i)-\mP_{n,\bp}(f)\pff\Mid\pffff\\
\nono\le & \EEE\poooo \pooo\sli_{f\in B_{n,\bp}}\Mid \sli_{i\in \NNN}p_i\poo f(Y_i)-\mP_{n,\bp}(f)\pff\Mid\pfff^2\pffff\\
\nono\le & (2N)^2\max_{f\in B_{n,\bp}}\;\EEE\poooo\Mid \sli_{i\in \NNN}p_i\poo f(Y_i)-\mP_{n,\bp}(f)\pff\Mid^2\pffff\\
\nono\le & (2N)^2 \max_{f\in B_{n,\bp}}\;\sli_{i\in \NNN}p_i^2\Var\poo f(Y_1)\pff\\ \le &(2NM \mmi \bp\mmi_2)^2.\label{haha}
\end{align}
Combining (\ref{ha}) and (\ref{haha}) yields, almost surely
$$\EEE\pooo \mmi G_n^M\mmi_\FF\Mid \mbeta_n\pfff\indic_{\SSS_n}(\mbeta_n)\le 2NM \mmi \mbeta_n\mmi_2+\e.$$
This concludes the proof, since $\mmi\mbeta_n\mmi_2\cvproba 0$ by assumption and since $\PPP\po \mbeta_n\notin \SSS_n\pf<\e$. $\Box$\lb
With Lemma \ref{ lem: GC sous un seuil arbitraire} at hand, the proof of Theorem \ref{T1} will be concluded as follows:
\begin{lem}\label{lem: deseuillage pour GC}
We have $$\lim_{M\rar \infty}\lsn d(G_n^M,G_n)=0.$$
\end{lem}
\textbf{Proof}: In view of (\ref{condition GC mbeta_n borne en proba}) it is sufficient to show that 
\beq \nono \forall \e>0,\;\mathop{\overline{\lim}}_{M\rar \infty}\lsn\PPP\pooo \EEE\poo F\indic_{\{F>M\}}\po Y_{1,n}\pf\Mid \mbeta_n\pff>\e\pfff\le \e.\eeq
This is immediate by (\ref{condition envelope GC}) combined with Markov's inequality.
$\Box$
\subsection{Proof of Theorem \ref{T2}}
By (\ref{ condition Donsker norme L2 egale a 1}) we can assume without loss of generality that $\mmi \mbeta_n\mmi_2\equiv 1$ for all $n$.
First note that (\ref{condition Donsker sur l'entropie crochet}) immediately implies
\beq \forall \dd>0, \pooo N_{[]}\po \dd,\FF,\norm_{\mP_{n,\mbeta_n},2}\pf\pfff_{n\geq 1} \text{is bounded in probability}\label{condition Donsker sur mfa(delta)}.\eeq
The proof of Theorem \ref{T2} follows the same directions as in that of Theorem 2 in \cite{Varron14Donsker}. The only crucial point that changes is that of proving the following asymptotic equicontinuity condition 
\beq  \lim_{\dd\rar 0} \; {\lsn}^{\PPP^*} \sup_{(f_1,f_2)\in \FF^2,\; \rho(f_1,f_2)< \dd}\mid G_n(f_1)-G_n(f_2)\mid =0,\label{equicontinuite cruciale}\eeq
which would be the only missing ingredient to complete the proof of Theorem 2.
Proving (\ref{equicontinuite cruciale}) will be achieved by conditioning upon $\mbeta_n$ and using the following chaining argument. It is an extension of usual chaining arguments for the bracketing entropy \cite[p. 286, Lemma 19.34]{VanderAsymptotic} to unbalanced empirical measures. Due to the fact that infinitely many weights are involved, only uniformly bounded classes of functions are treated here for simplicity. This will be largely sufficient for our purposes.
\begin{lem}\label{lem: chaining pour des empiriques desequilibrees} 
Let $\bp\in \ell^{1,+}$ such that $\mmi \bp\mmi_2=1$ let $Q$ be a probability measure and let $\GG$ be a uniformly bounded pointwise measurable class of functions with countable separant $\GG_0$. Let $\dd\in (0,\infty]$ be such that
\begin{align} 
&\sup_{g\in \GG}\mmi g\mmi_{Q,2}\le \dd\text{ and }\label{bep variance}\\ 
&\sup_{g\in \GG,\; y\in \mfX}\;\mid g(y)\mid\le \mmi \bp\mmi_{\infty}^{-1}\mfa(\dd,Q),\text{ where}\label{bep norme infinie}\\
\nono &\mfa(\dd,Q):=\dd/\sqrt{\log N_{[]}\po \dd,\GG,\norm_{Q,2}\pf}.
\end{align}
Then for any i.i.d sequence $(Y_i)_{i\in \NNN}$ with distribution $Q$, we have 
\beq \EEE\pooo \sup_{g\in \GG}\Mid \sli_{i\in \NNN}p_i\poo g(Y_i)-Q(g)\pff\Mid\pfff\le \mathfrak{C}_1
J_{[]}\po \dd,\GG,Q\pf,\label{ assertion du chaining unbalanced}\eeq
where $\mathfrak{C}_1$ is a universal constant. 
\end{lem}
\textbf{Proof}: 
We shall use the notations
\begin{align*}
\Delta^2(\GG,Q):=&\sup_{g\in \GG}Q(g^2),\;\text{and }\Gam(\GG):=\sup_{g\in \GG,\; y\in \mfX} \mid g(y)\mid.
\end{align*}
Given a finite class of functions $\tilde{\GG}$ and given $m\geq 1$ we have, by combining Lemmas 2.2.9 and 2.2.10 in \cite[p. 102]{Vander}:
\begin{align}
\nono &\EEE\poooo \max_{g\in \tilde{\GG}} \Mid \sli_{i=0}^m p_i \poo g(Y_i)-Q(g)\pff\Mid\pffff\\
\nono \le & 24 \coooo\sqrt{\sli_{i=0}^m p_i^2\Delta^2(\tilde{\GG},Q)\log(1+\sharp \tilde{\GG})}+ \max_{i\le m}\mid p_i\mid \Gam(\tilde{\GG})\log(1+\sharp \tilde{\GG})\cffff,\\
 \le &24 \coooo \mmi \bp\mmi_2\Delta(\tilde{\GG},Q)\sqrt{\log(1+\sharp \tilde{\GG})}+\mmi \bp\mmi_{\infty}\Gam(\tilde{\GG})\log(1+\sharp \tilde{\GG})\cffff,\label{papou}
  \end{align}
  where the possible choice of factor 24 was actually shown in \cite[p. 285, Lemma 19.33]{VanderAsymptotic}. Since (\ref{papou}) does not depend upon $m$ we then have, as soon as $\Gam(\tilde{\GG})<\infty$ (and recalling that $\mmi \bp\mmi_2=1$):
  \begin{align}
&\nono\EEE\poooo \max_{g\in \tilde{\GG}} \Mid \sli_{i\in \NNN} p_i \poo g(Y_i)-Q(g)\pff\Mid\pffff \\
 =&\lim_{m\rar \infty} \EEE\poooo \max_{g\in \tilde{\GG}} \Mid \sli_{i=0}^m p_i \poo g(Y_i)-Q(g)\pff\Mid\pffff\label{pa}\\ \le &24 \cooo \Delta(\tilde{\GG},Q)\sqrt{\log(1+\sharp \tilde{\GG})}+\mmi \bp\mmi_{\infty}\Gam(\tilde{\GG})\log(1+\sharp \tilde{\GG})\cfff, \label{inegalite maximale}
  \end{align}
 where $(\ref{pa})$ is an application of the dominated convergence theorem, since all the involved random variables are bounded by $2\Gam(\tilde{\GG})$. Now with $(\ref{inegalite maximale})$ at hand, the remainder of the proof is as follows: a careful look at all the arguments of the proof of \cite[p. 286, Lemma 19.34]{VanderAsymptotic} - noting that their truncating argument is not needed here - shows that the latter are still true with the systematic formal change of $\sqrt{n}$ by $\mmi \bp \mmi_{\infty}^{-1}$. $\Box$\lb
We can now start our proof of (\ref{equicontinuite cruciale}). First fix $\e>0$. Using (\ref{condition Donsker sur l'entropie crochet}) and (\ref{ condition sur la semimetrique rho}) there exist $\dd_1,\dd_2>0$ and $n_0$ such that for all $n\geq n_0$ we have $1-\e\le \PPP_*\po \mbeta_n\in \SSS'_n\pf$, where $\SSS'_n$ is the set of all $\bp\in \ell^{1,+}$ satisfying the following conditions:
\begin{align} 
&2\sqrt{2}\mathfrak{C}_1J_{[]}\poo \frac{\dd_1}{2},\FF,\mP_{n,\bp}\pff\le \e\label{ fil1}\\
\nono &\sup_{f\in \FF_{\dd_2}}\mmi f\mmi_{\mP_{n,\bp},2}<\dd_1,\text{ where  }\\
\nono &\FF_{\dd_2}:= \aoo f_1-f_2,\; (f_1,f_2)\in \FF^2,\; \rho(f_1,f_2)<\dd_2\aff,
\end{align}
and where $\mathfrak{C}_1$ denotes the universal constant in (\ref{ assertion du chaining unbalanced}).
Now fix $\bp$, write 
\begin{align*}
T(\bp):=& \mmi \bp\mmi_{\infty}^{-1}\mfa(\dd_1,\mP_{n,\bp})\indic_{\{\mmi \bp\mmi_{\infty}>0\}},\text{ and define}\\
\FF_{\bp,\dd_1}:= &\aoo (f_1-f_2)\indic_{\{F\le T(\bp)\}},\; (f_1,f_2)\in \FF^2,\; \mmi f_1-f_2\mmi_{\mP_{n,\bp},2}<\dd_1\aff.
\end{align*}
Next, apply Lemma \ref{lem: chaining pour des empiriques desequilibrees} for fixed $\bp \in \SSS'_n$ to obtain (noticing that $\FF_{\bp,\dd_1}$ satisfies (\ref{bep variance}) and (\ref{bep norme infinie}) for the choice of $Q:=\mP_{n,\bp}$ and $\dd:=\dd_1$)
\begin{align}
\nono \EEE\poooo \sup_{f\in \FF_{\bp,\dd_1}} \Mid \sli_{i\in \NNN} p_i\poo f (Y_i)-\mP_{n,\bp}\po f\pf\pff\Mid\pffff\le &\mathfrak{C}_1J_{[]}\poo \dd_1,\FF_{\bp,\dd_1},\mP_{n,\bp}\pff\,\\
\le &2\sqrt{2}\mathfrak{C}_1J_{[]}\poo \frac{\dd_1}{2},\FF,\mP_{n,\bp}\pff, \label{mpap}
\end{align}
where the $Y_i$ are i.i.d with law $\mP_{n,\bp}$ and where (\ref{mpap}) is a consequence of (\ref{ fil1}) and standard comparisons of entropy numbers. Now since the latter inequality is valid for all $\bp \in \SSS'_n$ we have 
\begin{align}
\EEE\pooo \mmi G_n^T\mmi_{\FF_{\mbeta_n,\dd_1}}\Mid \mbeta_n\pfff\indic_{\SSS'_n}\po \mbeta_n\pf\le \e\indic_{\SSS'_n}\po \mbeta_n\pf,\text{ almost surely.}\label{oh}
\end{align}
Note that the measurability $\mmi G_n^T\mmi_{\FF_{\mbeta_n,\dd_1}}$ is not immediate at all, but can be proved using the same arguments as in \cite[proof of Proposition 4.2]{Varron14Donsker}.
In view of (\ref{oh}), and since $\FF_{\dd_2}\subset \FF_{\bp,\dd_1}$ for $\bp\in \SSS'_n$, the proof of (\ref{equicontinuite cruciale}) will be completed if we prove the following lemma.
\begin{lem}
We have
$d(G_n^T,G_n)\rar 0$ as $\nif$.
\end{lem}

\textbf{Proof}: From (\ref{inegalite entre les esperance condtionnelles pour le reste de cauchy}) we have (noting that $\mmi \beta_n\mmi_2\equiv1$ implies $\mmi\beta_n\mmi_{\infty}>0$ a.s.)
\begin{align}
\nono &d\po G_n^{T},G_n\pf\\
\le &\EEE\poooo \arctan\pooo 2\mmi \mbeta_{n}\mmi_1\EEE\poo F\indic_{\{F>T(\mbeta_n)\}}(Y_{1,n})\mid \mbeta_n\pff\pfff\pffff \label{hao}\\
\nono =&\EEE\poooo \arctan\pooo 2\mmi \mbeta_{n}\mmi_1\EEE\poo F\indic_{\{F>T(\mbeta_n)\}}(Y_{1,n})\mid \mbeta_n\pff\pfff\pffff\\
\nono  \le & \EEE\poooo \arctan\pooo 2\frac{\mmi \mbeta_{n}\mmi_1}{T(\mbeta_n)^{p-1}}\EEE\poo F^p\indic_{\{F>T(\mbeta_n)\}}(Y_{1,n})\mid \mbeta_n\pff\pfff\pffff\\
\nono = & \EEE\poooo \arctan\pooo 2\frac{\mmi \mbeta_{n}\mmi_1\times \mmi \mbeta_n\mmi_{\infty}^{p-1}}{\mfa(\dd_1,\mP_{n,\mbeta_n})^{p-1}}\EEE\poo F^p\indic_{\{F>T(\mbeta_n)\}}(Y_{1,n})\mid \mbeta_n\pff\pfff\pffff.
\end{align}
Since, by (\ref{condition Donsker produit l1linfty borne}) and (\ref{condition Donsker sur mfa(delta)}), the sequence $\mmi \mbeta_{n}\mmi_1\times \po\mmi \mbeta_n\mmi_{\infty}/\mfa(\dd_1,\mP_{n,\mbeta_n})\pf^{p-1}$ is bounded in probability, it only remains to prove that 
\beq E_n:=\EEE\pooo F^p\indic_{\{F>T(\mbeta_n)\}}(Y_{1,n})\mid \mbeta_n\pfff\cvproba 0\label{gt}.\eeq
 To prove this, fix $\e>0$ and choose $M$ large enough so that $$\EEE\poo F^p\indic_{\{F>M\}}(Y_{1,n})\pff\le \e^2,$$
for all $n\geq 1$, which is possible by (\ref{ condition envelope Donsker}). Next apply Markov's inequality to $E_n$ on the set $\{ T(\mbeta_n)> M\}$ and then note that $\PPP\po T(\mbeta_n)\le M\pf\rar 0$ by (\ref{condition Donsker linfty tend vers 0}) and (\ref{condition Donsker sur mfa(delta)}). To conclude the proof, let us now consider the isolated case where $(\ref{condition Donsker produit l1linfty borne})$ and $(\ref{ condition envelope Donsker})$ are removed from the set of assumptions of Theorem \ref{T2}, but $\FF$ is uniformly bounded, i.e., $F\le M$ for some constant $M>0$. Then a look at (\ref{hao}) immediately yields the claim, noticing that $T(\mbeta_n)\cvproba \infty$. $\Box$
\subsection{Proof of Corollary \ref{ coro:posterior}}\label{subsection preuve posterior}With (\ref{representation de Sethuraman}) in mind, let us define 
\beq A:=\aoo (x_n)_{n\geq 1},\; \Mmi \alp_{\xxn}-\mP_0\Mmi_{\FF} \rar 0\aff.\label{definition de A}\eeq
The class of indicators of subsets of a countable set is universally Glivenko--Cantelli - see, e.g., \cite[p. 217, Remark 6.4.3]{Dudley}. Therefore, since $\tht_n\rar 0$, the triangle inequality entails $\mP_0^{\NNN^*}(A)=1$. Now take an arbitrary sequence $(x_n)_{n\geq 1}\in A$. We shall apply Theorem \ref{T1} to the sequence $Post_n\xxn$. In this setup we have $\mP_{n,\bp}=\mP_n=\alp_{\xxn}$ for all $\bp \in \ell^1$, and $$\beta_{i,n}:=V_{i,n}\proli_{j=0}^{i-1}(1-V_{j,n}),\; i\in \NNN,\; n\geq 1,$$
with $(V_{i,n})_{i\in \NNN}\leadsto Beta(1,M+n)^{\otimes \NNN}$. To prove (\ref{condition bracket GC}) let us first remark that if $\llbracket f^-,f^+\rrbracket$ is a bracket between two indicator functions fulfilling $\mmi f^+-f^-\mmi_{\mP_{\xxn},1}\le \e$ then $\mmi f^+-f^-\mmi_{\alp_{\xxn},1}\le \tht_n+(1-\tht_n)\e$. Moreover, since the pointwise supremum/infimum of a set of indicator functions is itself an indicator function, any covering of $\FF$ by brackets can be converted into another covering with the same number of brackets, each of one between two indicator functions. Hence, since $\tht_n\rar 0$, we conclude that it is sufficient to prove that $N_{[]}(\e,\FF,\mP_{\xxn})$ is a bounded sequence for fixed $\e>0$. This is done as follows: let us first choose a finite set $C_0\subset \mfX$ such that $\mP_0(C_0)>1-\e$. Then by definition of $A$ one has $\mP_{\xxn}(C_0)>1-\e$ for all large enough $n$.
\beq\forall C\subset \mfX, C\cap C_0\subset C\subset \po  C\cap C_0\pf \cup C_0^c\label{inclusion}.\eeq
Hence the the finite collection 
$$ \aoo \llbracket \mathds{1}_{C},\mathds{1}_{C\cup C_0^c}\rrbracket,\; C\subset C_0\aff$$
defines a covering of $2^{\sharp C_0}$ brackets having $\norm_{\mP_{\xxn},1}$ diameters less than $\e$. This proves that 
$N_{[]}(\e,\FF,\mP_{\xxn})\le 2^{\sharp C_0}$ for all large $n$, and hence proves (\ref{condition bracket GC}). Now 
conditions (\ref{condition envelope GC}) and  (\ref{condition GC mbeta_n borne en proba}) are immediate since $\FF$ is uniformly bounded and $\mmi \mbeta_n\mmi_1\equiv 1$ - see, e.g. \cite[p. 112]{HjortLivre}. Finally, standard calculus on beta distributions shows that 
$\EEE\po \mmi\mbeta_n\mmi_2^2\pf\sim n^{-1}$, from where one can apply Theorem \ref{T1} and conclude the proof.

\subsection{Proof of Corollary \ref{coro: posterior et bernstein}}\label{sous section preuve Bernstein}
We shall now assume without loss of generality that the support of $\mP_0$ is infinite.
\subsubsection{Two preliminary results}
Theorem \ref{theo Durst} states that the finiteness of $\Sigma_{y\in \mfX}\sqrt{\mP_0(\{y\})}$ is equivalent to that of $J_{[]}(\infty,\FF,\mP_0)$. Our next lemma goes one step further: it shows that it is possible to control the magnitude of $J_{[]}(\dd,\FF,\mP_0)$, for small $\dd>0$, by "tail" sums of the $\sqrt{\mP_0(\{y\})}$.
\begin{lem}\label{lem des sommes de queues}
Define, for $k\in \NNN:$ 
\beq \mathbf{j}_{\mP_0}(k):=\min\aoo J\in \NNN, \;\sli_{y\in \mfX\;:\; \pl\le 16^{-J} }\pl\le   4^{-k}\aff\label{definition de j(k)}.\eeq
Then, for all $p\geq 1$ we have, for a universal constant $\mathfrak{C}_2$
$$J_{[]}\po 2^{-(p-1)},\FF,\mP_0\pf\le \mathfrak{C}_2 \sqrt{\sli_{y\in  \mfX}\sqrt{\mP_0(\{y\})}}\times \sqrt{\sli_{y: \; \pl\le 16^{-\mathbf{j}_{\mP_0}(p)+1}}\sqrt{\mP_0(\{y\})}}.$$
Moreover if the support of $\mP_0$ is infinite we have $\mathbf{j}_{\mP_0}(p)\rar \infty$ as $p\rar \infty$.
\end{lem}
\textbf{Proof }: The very last statement is obvious. We shall now write $\mathbf{j}(\cdot)$ instead of $\mathbf{j}_{\mP_0}(\cdot)$ for concision. The proof consists in enriching the arguments of Dudley \cite[p. 245-246]{Dudley} with additional analytical precisions. We shall hence borrow his notations. First, for $j\in \NNN$ write 
\beq \nono A_j:=\aoo y\in  \mfX,\; 16^{-j-1}<\mP_0(\{y\})\le 16^{-j}\aff,\text{ and }r_j:=\sharp A_j\label{definition de Aj}.\eeq Now define the following maps on $\NNN$
\begin{align}
\nono&m(\cdot):\; k\rar \sli_{j=0}^{\mathbf{j}(k)}r_j=\sharp \bculi_{j=0}^{\mathbf{j}(k)}A_j,\\
\nono &k(\cdot):\; J\rar \min\aoo p\geq 1,\:4^{-p}< \sli_{y:\; \pl\le 16^{-J}}\pl\aff,\\
\nono &\kappa(\cdot):\; k\rar \min\aoo \kappa \in \NNN,\; \mathbf{j}(\kappa)=\mathbf{j}(k)\aff.
\end{align}
For consistency of notations in the following calculus, we shall also define $k(-1):=0$. Note that, writing $\KK$ for the range of $\kappa(\cdot)$, the map $\mathbf{j}(\cdot)$ is one to one on $\KK$.\lb
Fix $k\geq 1$. Similarly as in Dudley \cite[p. 245-246]{Dudley} we see that, for fixed $k\geq 1$, one can use the same arguments as for (\ref{inclusion}), with the formal replacement of $\e$ by $4^{-k}$ and $ C_0$ by 
$$C_k:=  \bculi_{j=0}^{\mathbf{j}(k)}A_j,$$ 
which satisfies $\mP_0(C_k)\geq 1-4^{-k}$ by (\ref{definition de j(k)}). This implies

\beq \forall k\geq 1,\;  N_{[]}\po 2^{-k},\FF,\norm_{\mP_0,2}\pf \le 2^{m(k)}.\label{borne crochet mk}\eeq 
Now, for any $p\geq 1$, by monotonicity of the involved functions:
\begin{align*}
J_{[]}\po 2^{-(p-1)},\FF,\mP_0\pf
\le& \sli_{k\geq p}\sqrt{\log N_{[]}\po 2^{-k},\FF,\norm_{\mP_0,2}\pf}\po 2^{-(k-1)}-2^{-k}\pf\\
\le& \sqrt{\log (2)}\sli_{k\geq p}\frac{\sqrt{m(k)}}{2^k} \text{ by }(\ref{borne crochet mk}).
\end{align*}
Next, fix $p\geq 1$ and write\begin{align*}
\sli_{k\geq p}\frac{\sqrt{m(k)}}{2^k}\le& \sli_{k\geq p}\sli_{j=0}^{\mathbf{j}(k)}\sqrt{r_j}\;2^{-k}\\
\le &\sli_{k\geq p}\sli_{j=0}^{\mathbf{j}(k)}2\sqrt{\sli_{y\in  A_j}\sqrt{\pl}}\;2^{j-k}\text{, since }r_j\le \sum_{y\in A_j}4^{j+1}\sqrt{\pl}\\
=&2\sli_{j\geq 0}\sqrt{\sli_{y\in  A_j}\sqrt{\pl}}\mathop{\sli_{k:\;k\geq p,}}_{\mathbf{j}(k)\geq j}2^{j-k}\\
 \le& 2\sqrt{\sli_{j\geq 0}\sli_{y\in  A_j}\sqrt{\pl}}\times\sqrt{\sli_{j\geq 0}\po \mathop{\sli_{k:\;k\geq p,}}_{\mathbf{j}(k)\geq j}2^{j-k}\pf^2},\text{ using Cauchy--Schwartz}\\
 =& 2\sqrt{\sli_{y\in  \mfX}\sqrt{\pl}}\times \sqrt{\sli_{j\geq 0}\po \mathop{\sli_{k\geq p,}}_{\mathbf{j}(k)\geq j}2^{j-k}\pf^2}.
\end{align*}
Now we have
\begin{align*}
&\sli_{j\geq 0}\po \mathop{\sli_{k:\;k\geq p,}}_{\mathbf{j}(k)\geq j}2^{j-k}\pf^2\\
\le &4\sli_{j\geq 0}4^{j-k(j-1)}\wedge 4^{j-p},\text{ since } \mathbf{j}(k)\geq j\text{ implies } k\geq k(j-1)\\
=& 4\sli_{k\geq 0}\;\sli_{j: \;\mathbf{j}(k-1)\le j-1<\mathbf{j}(k)}4^{j-k}\wedge 4^{j-p} \text{, since }k(j)=k \text{ for }\mathbf{j}(k-1)\le j< \mathbf{j}(k)\\
=&4\sli_{k\le p} \;\sli_{j:\;\mathbf{j}(k-1)\le j-1< \mathbf{j}(k)}4^{j-p}+4\sli_{k\geq p+1} \;\sli_{j:\;\mathbf{j}(k-1)\le j-1< \mathbf{j}(k)}4^{j-k}\\
\le &8\coo4^{\mathbf{j}(p)-p}+\sli_{k\geq p+1}4^{\mathbf{j}(k)-k}\cff\\
=& 8\sli_{k\geq p}4^{\mathbf{j}(k)-k}\\
\le &8\sli_{k\geq p}4^{\mathbf{j}(k)-k+\kappa(k)}\times \sli_{\ell\;:\;\pl\le 16^{-\mathbf{j}(k)+1}}\pl,\text{ by }(\ref{definition de j(k)})\text{ and since  }\mathbf{j}(k)=\mathbf{j}(\kappa(k))\\
= &8\sli_{k\geq p}4^{\mathbf{j}(k)-k+\kappa(k)}\times \sli_{j\geq \mathbf{j}(k)-1}\;\sli_{\ell\;\in A_j}\pl\\
\le & 8\sli_{j\geq 0}\poo \sli_{y\in  A_j}\pl\pff\poo\sli_{k:\;k\geq p,\; \mathbf{j}(k)-1\le j} 4^{\mathbf{j}(k)-k+\kappa(k)}\pff.
\end{align*}
Now notice that, when $\mathbf{j}(p)>j+1$ the set of indices $\{k\geq p,\; \mathbf{j}(k)-1\le j\}$ is empty, from where\begin{align*}
&\sli_{j\geq 0} \poo\sli_{y\in  A_j}\pl\pff \poo \sli_{k:\;k\geq p,\; \mathbf{j}(k)-1\le j} 4^{\mathbf{j}(k)-k+\kappa(k)}\pff\\ 
\le & \sli_{j\geq \mathbf{j}(p)-1}\poo \sli_{y\in  A_j}\pl\pff \poo \sli_{k:\; \mathbf{j}(k)-1\le j} 4^{\mathbf{j}(k)+\kappa(k)-k}\pff\\
\le& \sli_{j\geq \mathbf{j}(p)-1} \poo \sli_{y\in  A_j}\pl \pff \poo \sli_{k'\in \KK,\; \mathbf{j}(k')\le j+1}4^{\mathbf{j}(k')+k'}\sli_{k:\; \kappa(k)=k'}4^{-k}\pff\\
\le&4\sli_{j\geq \mathbf{j}(p)-1} \poo \sli_{y\in  A_j}\pl\pff \poo\sli_{k'\in \KK,\; \mathbf{j}(k')\le j+1}4^{\mathbf{j}(k')}\pff\text{, since }\kappa(k)=k'\text{ implies } k\geq k'\\
\le &32\sli_{j\geq \mathbf{j}(p)-1} \sli_{y\in  A_j}\pl 4^j,\text{ since }\kappa(\cdot)\text{ is one to one on }\KK\\
\le &32 \sli_{j\geq \mathbf{j}(p)-1}\sli_{y\in  A_j}\sqrt{\pl},\text{ since }y\in  A_j\text{ implies }\pl 4^j\le \sqrt{\pl}\\
=&32\sli_{y: \; \pl\le 16^{-\mathbf{j}(p)+1}}\sqrt{\pl}.\end{align*}
This concludes the proof.$\Box$\lb
Our second preliminary result is as follows.
\begin{lem}\label{lem application du multiplier CLT}
Write 
$$I_{\e}:=\ao y\in  \mfX,\; \pl \le  \e\af,\e\in \QQQ^+.$$
Then for $\mP_0^{\otimes \NNN^*}$-almost any sequence $(x_n)_{n\geq 1}$ we have:
\beq \forall \e\in \QQQ^+,\; \limn \sli_{y\in  I_\e}\sqrt{\mP_{\xxn}(\{y\})}= \sli_{y\in   I_\e}\sqrt{\mP_0(\{y\})}.\eeq
\end{lem}
\textbf{Proof:} Since the class $\FF$ is $\mP_0$-Donsker and admits a square integrable envelope ($F\equiv1$), the conditional multiplier Donsker theorem applies for a suitable i.i.d. standard normal sequence $(\xi_n)_{n\geq 1}$ - see, e.g., \cite[p. 183, Theorem 2.9.7]{Vander}. Hence for $\mP_0^{\NNN^*}$-almost every sequence $(x_n)_{n\geq 1}$ we have - recalling that $\WW_{\mP_0}$ stands for the $L^2(\mP_0)$-isonormal Gaussian process indexed by $\FF$:
\begin{align}
\pooo \WW_{\xxn}(f)\pfff_{f\in \FF} \cvloi \poo \WW_{\mP_0}(f)\pff_{f\in \FF},\label{conditional CLT}\text{ where }\\
\nono\WW_{\xxn}(f):=\frac{1}{\sqrt{n}}\sliin \xi_if(x_i),\; f\in \FF.
\end{align}
Here the weak convergence holds in the sense of Hoffman-J\H{o}rgensen holds taking the underlying probability space as the canonical product space for $(\xi_n)_{n\geq 1}$. Moreover the involved processes are Gaussian, hence weak convergence implies convergence of first moments of absolute suprema. As a consequence, for such a sequence $(x_n)_{n\geq 1}$ fulfilling (\ref{conditional CLT}) we have, for all $\e\in \QQQ^+$
$$\EEE\poooo2\sup_{A\subset I_\e }\frac{1}{\sqrt{n}}\Mid \sliin \xi_i\mathds{1}_{A}(x_i)\Mid- \frac{1}{\sqrt{n}}\Mid \sliin \xi_i\mathds{1}_{I_\e}(x_i)\Mid\pffff\rar \EEE\poooo2\sup_{A\subset I_\e } \Mid \WW_{\mP_0}\po \mathds{1}_A\pf\Mid-\Mid \WW_{\mP_0}\po\mathds{1}_{I_\e}\pf\Mid\pffff.$$
Finally, by the standard equality 
$$\sup_{A\subset I_{\e}}\Mid \sli_{y\in  A}g(y)\Mid =\frac{1}{2}\poo \sli_{y\in  I_{\e}}\mid g(y)\mid+\Mid \sli_{y\in  I_{\e}}g(y)\Mid\pff,$$
we have (with $g(y):=n^{-1/2}\sliin \xi_i\mathds{1}_{\{y\}}(x_i)$)
\begin{align*} 
&\EEE\poooo2\sup_{A\subset I_\e }\frac{1}{\sqrt{n}}\Mid \sliin \xi_i\mathds{1}_{A}(x_i)\Mid- \frac{1}{\sqrt{n}}\Mid \sliin \xi_i\mathds{1}_{I_\e}(x_i)\Mid\pffff\\
=&
\EEE\poooo \sli_{y\in  I_\e}\Mid \frac{1}{\sqrt{n}} \sliin \xi_i\mathds{1}_{\{y\}}(x_i)\Mid
\pffff\\
=&\sli_{y\in  I_{\e}}\sqrt{\mP_{\xxn}(\{y\})},
\end{align*}
and similarly (with now $g(y):=\WW_{\mP_0}(\mathds{1}_{\{y\}})$)
\begin{align*}
\EEE\pooo 2\sup_{A\subset I_\e } \Mid \WW_{\mP_0}\po \mathds{1}_A\pf\Mid-\Mid \WW_{\mP_0}\po\mathds{1}_{I_\e}\pf\Mid\pfff=\sli_{y\in  I_\e}\sqrt{\mP_0(\{y\})}\,
\end{align*}
which concludes the proof.$\Box$
\subsubsection{Use of Theorem \ref{T2} }
Recall that $A$ was defined in (\ref{definition de A}) and has probability one.
Let us consider the set 
\begin{align*}
B:=& \aooo (x_n)_{n\geq 1},\; \forall \e\in \QQQ^+,\; \lsn \sli_{y:\; \alp_{\xxn}(\{y\})\le \e}\sqrt{\alp_{\xxn}(\{y\})}\le  \sli_{y\in  I_{2\e}}\sqrt{\mP_0(\{y\})}\afff.
\end{align*}
We have, for any $n\geq 1$ and $\e\in \QQQ^+$ (using $\sqrt{a+b}\le \sqrt{a}+\sqrt{b}$)
\begin{align*}
&\sli_{y:\; \alp_{\xxn}(\{y\})\le \e}\sqrt{\alp_{\xxn}(\{y\})}\\
\le& \sqrt{\tht_n}\sli_{y:\; \alp_{\xxn}(\{y\})\le \e}\sqrt{\alp(\{y\})}+\sqrt{1-\tht_n}\sli_{y:\; \alp_{\xxn}(\{y\})\le \e}\sqrt{\mP_{\xxn}(\{y\})}.
\end{align*}
Now if $(x_n)_{n\geq 1}$ belongs to $A$ and since $\FF$ induces the total variation distance we have, for all $n$ large enough
$$\ao y\in  \mfX,\; \alp_{\xxn}(\{y\})\le \e\af\subset \ao y\in  \mfX,\; \mP_0(\{y\})\le 2\e\af=I_{2\e}.$$
We hence conclude that $\mP_0^{\NNN^*}(A\cap B)=1$ - recalling that $\tht_n\rar 0$ and $(DDB(\alp))$ holds. Let us now consider a sequence $(x_n)_{n\geq 1}\in A\cap B$. Similarly as in \S \ref{subsection preuve posterior}, we shall  prove Corollary \ref{coro: posterior et bernstein} by verifying all the assumptions of Theorem \ref{T2}, for the choice of  $\mPnp=\mP_n:=\mP_{\xxn}$.
Because the class $\FF$ is uniformly bounded by $1$, the conditions upon 
$$\beta_{i,n}:=\sqrt{n}V_{i,n}\proli_{j=0}^{i-1}(1-V_{j,n}),\; i\in \NNN,$$
 that we need to check are (\ref{ condition Donsker norme L2 egale a 1}) and $(\ref{condition Donsker linfty tend vers 0})$, or equivalently 
$$\mmi \mbeta_n\mmi_2\cvproba 1,\text{ and }\mmi \mbeta_n\mmi_4\cvproba 0.$$
These are respectively proved by direct computations of expectations and variances. It now remains to verify (\ref{condition Donsker sur l'entropie crochet}) and (\ref{ condition sur la semimetrique rho}). 
By definition of $A$ and since $\tht_n\rar0$, the sequence $\mP_n$ obviously fulfills (\ref{SW}) and therefore satisfies (\ref{ condition sur la semimetrique rho}). Now in view of Lemma \ref{lem des sommes de queues}, assertion (\ref{condition Donsker sur l'entropie crochet}) will be proved if we show that 
\beq\lim_{p\rar \infty}\; \lsn \sli_{y: \;\alp_{\xxn}(\{y\})\le 16^{-\mathbf{j}_{\mP_n}(p)+1}}\sqrt{\mP_n(\{y\})}=0.\label{a prouver}\eeq
\begin{lem}\label{lem1}
Take $(x_n)_{n\geq 1}\in A$. For any $p\geq 1$ we have $\mathbf{j}_{\mP_n}(p)\geq \mathbf{j}_{\mP_0}(p)$ for all large enough $n$.
\end{lem}
\textbf{Proof}: Fix $p\geq 1$. By definition of $\mathbf{j}_{\mP_0}$ we have 
$$\sli_{y:\; \mP_0(\{y\})\le 16^{-\mathbf{j}_{\mP_0}(p)+1}}\mP_0(\{y\})> 4^{-p}.
$$
Now since $(x_n)_{n\geq 1}\in A$, we have, for all $n$ large enough:
$$\sli_{y:\; \mP_0(\{y\})\le 16^{-\mathbf{j}_{\mP_0}(p)+1}}\mP_n(\{y\})> 4^{-p},
$$
whence $\mathbf{j}_{\mP_0}(p)-1\le \mathbf{j}_{\mP_n}(p)-1$ by definition of $\mathbf{j}_{\mP_n}$. $\Box$\lb
Now applying Lemma \ref{lem1} we have, writing $\e(p):=16^{-\mathbf{j}_{\mP_0}(p)+1}$:
\begin{align*}
&\lim_{p\rar \infty} \lsn\;\; \sli_{y:\; \;\alp_{\xxn}(\{y\})\le 16^{-\mathbf{j}_{\mP_n}(p)+1}}\sqrt{\mP_n(\{y\})}\\
\le &\lim_{p\rar \infty} \lsn \;\;\sli_{y:\; \;\alp_{\xxn}(\{y\})\le \e(p)}\sqrt{\mP_n(\{y\})}\\
\le &\lim_{p\rar \infty} \lsn \;\;\sli_{y\in  I_{2\e(p)}}\sqrt{\mP_0(\{y\})},\text{ since }(x_n)_{n\geq 1}\in B\\
=&0,
\end{align*}
by $(DDB(\mP_0))$ together with $\lim \mathbf{j}_{\mP_0}(p)\rar \infty$. This proves (\ref{a prouver}) and we can now apply Theorem \ref{T2} to obtain
$$d_{BL}\poooo\sqrt{n}\poo Post_n\xxn-\alp_{\xxn}\pff,\GGG_{\alp_{\xxn}}\pffff\rar 0.$$
But since $(x_n)_{n\geq 1}\in A$ the sequence $\mP_n:=\alp_{\xxn}$ satisfies (\ref{SW}) from where (see \cite[Remark 2.2]{Varron14Donsker}):
$$\GGG_{\alp_{\xxn}}\cvloi \GGG_{\mP_0}\text{, in }\ell^{\infty}(\FF),$$
which concludes the proof of Corollary \ref{coro: posterior et bernstein}. $\Box$
 \subsection{Proof of Corollary \ref{coro Donsker local}}\label{sous section preuve de Donsker local}
Recall that $\mP_0$ denotes here the uniform distribution on $S$ and that $\mP_n$ has been defined in (\ref{definition Pn local}). Let $\VV$ be a neighborhood of $z$ on which $Z_1$ admits the density $\mf$. Since $S$ is bounded and $h_n\rar 0$ we have $z+h_nS\subset \VV$ for $n$ large enough. We may assume without loss of generality that this is the case for all $n\geq 1$.
\begin{lem}\label{lem convergence du PE randomise pour le PE local}
We have (taking here the convention $0/0=0$)
$$\poooo\sliin \frac{b_{i,n}}{\sqrt{\sliin b_{i,n}^2}}\poo f(Y_{i,n})-\mP_n(f)\pff\pffff_{f\in \FF}\cvloi \GGG_{\mP_0},$$
where $\GGG_{\mP_0}$ denotes the $\mP_0$ Brownian bridge.
\end{lem}
\textbf{Proof}:
Write $$\beta_{i,n}:=\frac{b_{i,n}}{\sqrt{\sliin b_{i,n}^2}},\text{ for }i=1,\ldots,n. $$ Since $\mf$ is continuous at $z$ we have 
\begin{align}
& a_n\sim \lab(S)\mf(z)h_n^d\label{an equivalent de f(z)hn^dlab(S)},\text{ from where }na_n\rar \infty\text{ and }a_n\rar 0.
\end{align}
This property ensures that the sequence $\mbeta_n$ satisfies (\ref{ condition Donsker norme L2 egale a 1}), (\ref{condition Donsker linfty tend vers 0}) and (\ref{condition Donsker produit l1linfty borne}) of Theorem 2 - taking $p:=2$ and recalling that $b_{i,n}\equiv b_{i,n}^2$. In order to verify (\ref{condition Donsker sur l'entropie crochet}) and (\ref{condition Donsker sur mfa(delta)}) we will now prove (\ref{condition de domination}), noting here that $\mP_{n,\bp}:=\mP_n$ for all $\bp\in \ell^{1,+}$. The usual change of variable $u=h_n^{-1}(v-z)$ in the next integrals gives, for an arbitrary non negative function $g$ with support included in $S$
\begin{align}
\nono \mP_n(g)=& \frac{1}{a_n}\ili_{z+h_nS}g\poo h_n^{-1}(v-z)\pff \mf(v) dv\\
\nono =&\frac{h_n^d}{a_n}\ili_{S}g(u)\mf(z+h_nu)du\\
\nono \le & \sup_{u\in S}\mf(z+h_n u)\;\frac{h_n^d}{a_n}\ili_{S} g(u)du\\
 \nono = & \sup_{u\in S}\mf(z+h_n u)\;\frac{h_n^d}{a_n}\lab(S)\mP_0(g).
\end{align}
This proves (\ref{condition de domination}) by applying that inequality to elements of the form $(f_1-f_2)^2$, $(f_1,f_2)\in \FF^2$ and recalling (\ref{an equivalent de f(z)hn^dlab(S)}) together with the continuity of $\mf$ at $z$. This also proves (\ref{ condition envelope Donsker}), taking $g:=F^2\indic_{\{F> M\}}$. Let us now verify (\ref{ condition sur la semimetrique rho}) by proving (\ref{SW}).
Using a calculus similar as above we have, for an arbitrary function $g\prec  (2F)^2\vee (2F)$ \begin{align}
&\nono \mf(z)\lab(S)\Mid \mP_n(g)-\mP_0(g)\Mid\\
\nono  =& \Mid \frac{\mf(z)\lab(S)h_n^d}{a_n} \ili_S g(u)\mf(z+h_nu) du-\ili_{S} g(u)\mf(z) du\Mid\\
\nono \le&\frac{\mf(z)\lab(S)h_n^d}{a_n}\times  \Mid  \ili_S g(u)\mf(z+h_nu) du-\mf(z)\ili_{S} g(u)du \Mid\\
\nono &\;\; +\Mid \frac{\mf(z)\lab(S)h_n^d}{a_n}-1\Mid \times \mf(z)\ili_{S} \mid g (u)\mid du \\
\nono \le & \frac{\mf(z)\lab(S)h_n^d}{a_n}\times \sup_{v\in S} \mid \mf(z+h_nv)-\mf(z)\mid\times \ili_S (2F)^2\vee (2F) du\\
&\;\; +  \Mid \frac{\mf(z)\lab(S)h_n^d}{a_n}-1\Mid \times \mf(z)\ili_{S}  (2F)^2\vee (2F) du\label{ borne 2},
\end{align}
which tends to zero independently of $g\prec (2F)^2\vee (2F)$. This proves $(\ref{SW})$ and concludes the proof of Lemma \ref{lem convergence du PE randomise pour le PE local}.$\Box$ \lb
Let us now continue the proof of Corollary \ref{coro Donsker local}. First, note that 
we have 
\begin{align}
&\sliin b_{i,n}^2\sim\lab(S)\mf(z)nh_n^d \text{ in probability, from where}\\
& \poooo\frac{1}{\sqrt{\mf(z)\lab(S)nh_n^d}}\sliin \beta_{i,n}\poo f(Y_{i,n})-\EEE\po f(Y_{i,n})\pf\pff\pffff_{f\in \FF}\cvloi \GGG_{\mP_0},
 \end{align}
 and hence that sequence of processes is asymptotically tight (see, e.g., \cite[p. 20, Definition 1.3.7]{Vander}). Now elementary probability calculus shows that 
\beq\frac{\sliin (b_{i,n}-a_n)}{\sqrt{\mf(z)\lab(S)nh_n^d}}\cvloi Z,\label{normalite asymptotique des bin}\eeq
where $Z$ is standard normal. Moreover, since $\mP_n$ satisfies (\ref{SW}) and since 
 $f\rar \mP_0(f)$ is continuous with respect to $\norm_{\mP_0,2}$, which makes $\FF$ totally bounded, the (deterministic) sequence $\mP_n(\cdot)$ is relatively compact in $\ell^{\infty}(\FF)$. This, combined with (\ref{normalite asymptotique des bin}), implies that  the sequence $R_n(\cdot)$ - defined in (\ref{definition de R_n(f)}) - is asymptotically tight, and hence so is $T_{n,h_n}(\cdot)$ by summation. It will hence be proved to converge to $\WW_{\mP_0}$ if we prove finite marginal convergences. This is done by elementary analysis of characteristic functions, using the change of variable $u=h_n^{-1}(v-z)$ in the integrals. We omit details. $\Box$

\section{Appendix: a minor proof}\label{sous section preuve de la mesurabilite}
In this section we prove the measurability properties claimed in \S \ref{section: results}.
\begin{lem}\label{lem de mesurabilite de brackets}
For fixed $r\geq 1$ and $n\geq 1$, the map $(\e,\bp)\rar N_{[]}\poo \e,\FF,\norm_{\mP_{n,\bp},r}\pff$ is Borel from $]0,\infty[\times \ell^1$ to $\RRR^+$.
As a consequence, the maps 
$$\bp \rar J_{[]}\po \dd,\FF,\mP_{n,\bp}\pf,\; \dd>0,$$
are Borel.
\end{lem}
\textbf{Proof}: Fix $r\geq 1$ and $n\geq 1$. Any bracket is closed for the the pointwise topology, i.e., the topology spanned by the evaluation maps $\ao \{f\rar f(y)\},\; y\in \mfX\af$. Hence so is any finite union of brackets that covers $\FF_0$. Since $\FF$ is included in the closure of $\FF_0$ for the pointwise topology, we deduce that
\beq \forall (\e,\bp)\in ]0,\infty[\times \ell^1,\; N_{[]}\poo\e,\FF,\norm_{\mP_{n,\bp},r}\pff=N_{[]}\poo\e,\FF_0,\norm_{\mP_{n,\bp},r}\pff.\label{egalite bracket FF et FF0}\eeq
Now the proof of Lemma \ref{lem de mesurabilite de brackets} boils down to proving the measurability of 
$$H:\; (\e,\bp)\rar N_{[]}\poo\e,\FF_0,\norm_{\mP_{n,\bp},r}\pff.$$
This is done by noting that, for any $K\in \NNN$, the set 
\beq B_K:=\aoo (f_j^-,f_j^+)_{j=1,\ldots,K}\in  {(\FF_0^2)}^K,\; \FF_0\subset \bculi_{j=1}^K \llbracket f_j^-,f_j^+\rrbracket\aff\nono \eeq
is countable, and that 
\begin{align*}
H(\e,\bp)>K\;\;\Leftrightarrow\;\; \forall (f_j^-,f_j^+)_{j=1,\ldots,K}\in  B_K,\; \exists j\in \{ 1,\ldots,K\},\; \mmi f_j^+-f_j^-\mmi_{\mP_{n,\bp},r}>\e,
\end{align*}
which yields the claimed result, since for fixed Borel non negative $g$, the map $\bp\rar \mmi g\mmi_{\mP_{n,\bp},r}$ is Borel (recall that $\{\mP_{n,\bp},\;\bp\in \ell^1\}$ is regular).$\Box$

\bibliographystyle{plain}
\bibliography{biblioAM-new,biblioNZ-new}
\end{document}